\pdfoutput=1

\documentclass[12pt]{amsart}
\usepackage{graphicx,amssymb,amscd,amsfonts,pinlabel,color,bm,times,overpic,psfrag}

\theoremstyle{plain}
\newtheorem{thm}{Theorem}[section]

\newtheorem{lemma}[thm]{Lemma}
\newtheorem{cor}[thm]{Corollary}

\theoremstyle{definition}

\theoremstyle{remark}

\numberwithin{equation}{section}

\title{On the Augmentation Categories of Positive Braid Closures}
\author{Michael Menke}

\begin{document}

\begin{abstract}
	In this paper we show that 0-resolution of a crossing in the Legendrian closure of a positive braid induces a cohomologically faithful $A_\infty$ functor on augmentation categories. In particular, we compute the bilinearized Legendrian contact cohomology of these knots for augmentations induced by 0-resolution.
\end{abstract}

\maketitle

\section{Introduction}
\label{Sec:intr}
\noindent
The $A_\infty$ augmentation category $Aug_-(\Lambda)$ of a Legendrian link $\Lambda \subset \mathbb{R}^3_{xyz}$ was first introduced in \cite{BC}. The objects of this category are augmentations of the link and the morphisms are chain complexes whose cohomology is known as bilinearized Legendrian contact cohomology. Bilinearized Legendrian contact cohomology can be used to distinguish exact Lagrangian fillings of a Legendrian knot up to Hamiltonian isomorphism and it provides an obstruction to Lagrangian concordances \cite{C}. In \cite{NRSSZ} a unital $A_\infty$ augmentation category $Aug_+(\Lambda)$ was constructed and used to show that augmentations of a Legendrian knot correspond to constructible sheaves on $\mathbb{R}^2_{xz}$ whose singular support at infinity lies on $\Lambda$. 
\\
\\
In this paper we will focus on the class of knots which are the closures of positive braids. Positive braid closures are a particular well-behaved class of knots whose Legendrian contact homology has been studied extensively in \cite{K}. Our goal is to study the effects of 0-resolution on the augmentation categories of positive braid closures. We will show that Lagrangian cobordisms corresponding to 0-resolution give rise to cohomologically faithful functors on augmentation categories. 
\\
\\
We will use $\mathbb{Z}/2$ coefficients throughout. Our main result is the following:

\begin{thm}
	Let $\Lambda_+ \subset \mathbb{R}^3$ be the Legendrian closure of a positive braid and let $\Lambda_-$ be the link obtained after resolving a crossing of $\Lambda_+$. Then there are cohomologically faithful $A_\infty$ \ functors $F^+: Aug_+(\Lambda_-) \to Aug_+(\Lambda_+)$ and $F^-: Aug_-(\Lambda_-) \to Aug_-(\Lambda_+)$ which are induced by 0-resolution. 
\end{thm}
\noindent
Using the $Aug_-$ version of Theorem 1.1 we obtain the following result for the bilinearized Legendrian contact cohomology of positive braids:

\begin{cor}
	Let $F^-$ be the functor from Theorem 1.1. Let $\epsilon_1$ and $\epsilon_2$ be $\mathbb{Z}/2$ graded augmentations of $\Lambda_-$. Then
\[LCH^*_{F^-(\epsilon_1),F^-(\epsilon_2)}(\Lambda_+) \simeq LCH^*_{\epsilon_1,\epsilon_2}(\Lambda_-) \oplus \mathbb{Z}/2[0] \]
where $\mathbb{Z}/2[0]$ denotes a copy of $\mathbb{Z}/2$ in degree 0.
\end{cor}

\section{Preliminaries}
\label{Sec:back}

\subsection{Positive Braids and Legendrian Contact Homology} 
\begin{figure}
	\centering
	\includegraphics[width=4in]{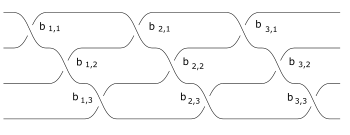}
	\caption{The (4,3) Torus Braid}
	\label{fig:Torus2}
\end{figure}

\begin{figure}
	\centering
	\includegraphics[width=4in]{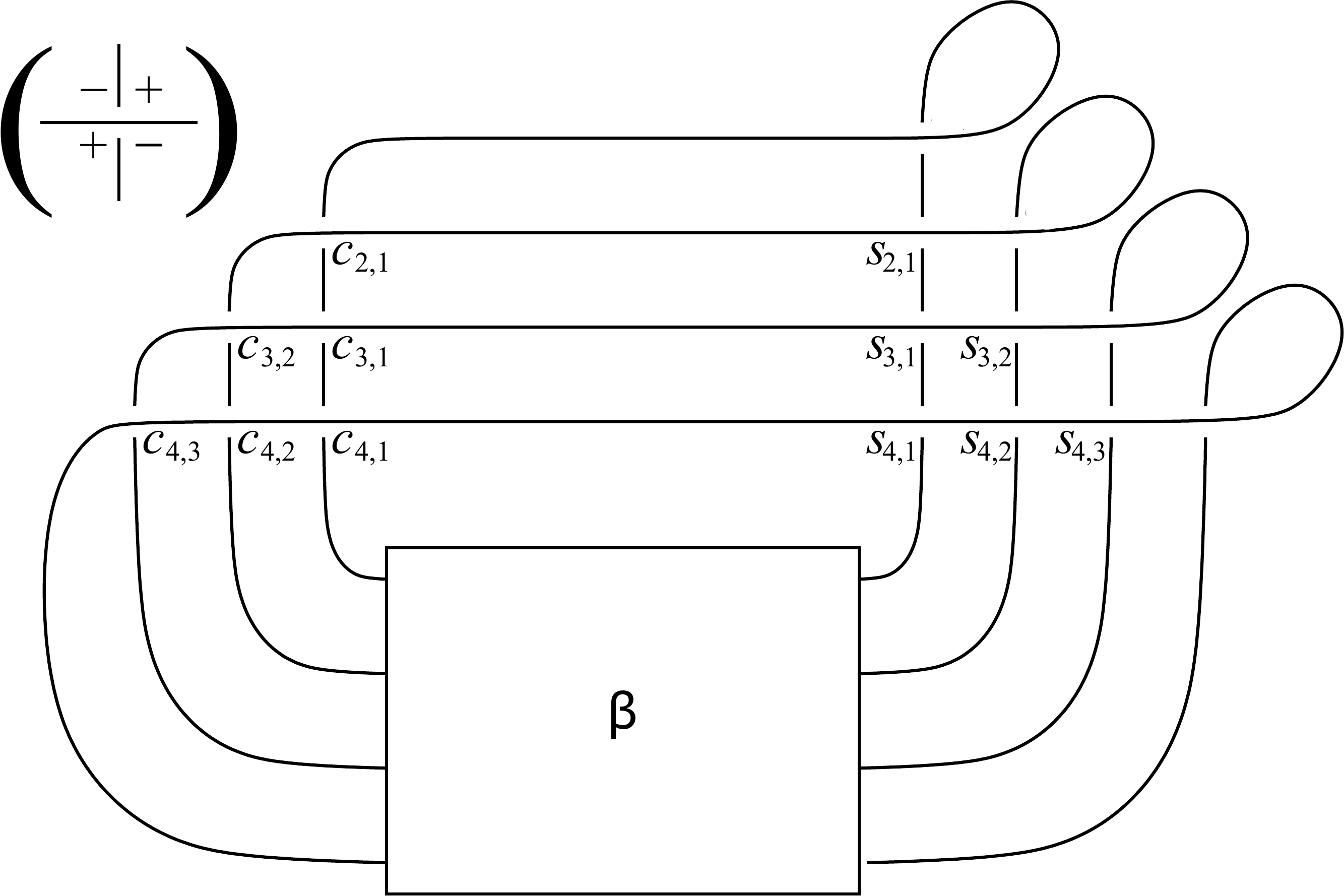}
	\caption{Lagrangian diagram of the Legendrian closure of the positive braid. The positive and negative corners of a crossing are labeled at the top left.}
	\label{fig:plat}
\end{figure}

A positive braid $\beta$ is a braid whose braid word consists entirely of right-handed twists. An example of a positive braid is the $(p,q)$ torus braid which can be represented by a diagram with $p$ strands and whose braid word is $(\sigma_1\sigma_2...\sigma_{p-1})^q$. The closure of a positive braid is obtained by attaching arcs which connect the beginning of strand $i$ to the right end of the strand which leaves the braid at position $i$. The closure of any positive braid has a Legendrian representative whose Lagrangian projection is shown in Figure \ref{fig:plat}. We will refer to this representative as the Legendrian representative.
\\
\\
We label the crossings in the braid portion of a torus link as in Figure \ref{fig:Torus2}. The crossings are $b_{i,j}$ where $1 \leq i \leq q$ and $1 \leq j \leq p-1$. Since every positive braid can be obtained as a sequence of 0-resolutions from a torus braid, every positive braid will inherit a labeling on its crossings from the corresponding torus knot.
\\
\\
To a Lagrangian diagram $\Lambda$ of a Legendrian knot we can associate a differential graded algebra $\mathcal{A}(\Lambda)$ with $\mathbb{Z}/2$ coefficients known as the Chekanov-Eliashberg DGA [Ch]. An alternative introduction can be found in [Et]. Let $\Lambda \subset \mathbb{R}^3$ where $\mathbb{R}^3$ has coordinates $x,y,z$. The projection of $\Lambda$ to the $xy$-plane is known as the Lagrangian projection. The generators of this DGA are the crossings of $\Lambda$ when viewed in the Lagrangian projection. These crossings correspond to Reeb chords of the knot. The differential of a crossing $a$ counts immersed disks in the $xy$-plane with punctures $a_1,b_1,\dots,b_n$ oriented counter-clockwise such that the $a_1$ puncture is at a positive convex corner of the crossing $a$, and the punctures $b_1,\dots,b_n$ are at negative convex corners of crossings in the diagram. The corners of a crossing are labeled as in Figure \ref{fig:plat}. The homology of this DGA is known as Legendrian contact homology. In the Chekanov-Eliashberg DGA of a positive braid the crossings $c_{i,j},b_{i,j}$ have degree 0 and the crossings $s_{i,j}$ have degree 1. More details on positive braid closures and their Legendrian contact homology can be found in \cite{K}.
\\
\\
A $\mathbb{Z}/2$-graded augmentation of a DGA $(A,\partial)$ is a map $\epsilon: A \to \mathbb{Z}/2$ such that $\epsilon \circ \partial = 0$ and $\epsilon(a) = 0$ unless $deg(a) = 0$. Throughout this paper by augmentation we will always mean a $\mathbb{Z}/2$-graded augmentation.

\subsection{0-resolution.}
In [EHK, Prop 6.17] it is shown that the Lagrangian cobordism corresponding to 0-resolution of a crossing $a$ of a link $\Lambda_+$ gives rise to a DGA map from $\mathcal{A}(\Lambda_{+})$ to $\mathcal{A}(\Lambda_{-})$, where $\Lambda_{+} = \Lambda_{-} \cup \{a\}$. This map is obtained by counting holomorphic disks with boundary on $\Lambda_+$ and two positive punctures, one of which is at the crossing $a$.
\\
\\
Let $d \neq a$ and denote by $M(d,a^k;\textbf{b})$ the moduli space of holomorphic disks in $\mathbb{R} \times \mathbb{R}^3$ and boundary on $\mathbb{R} \times \Lambda_+$ modulo biholomorpism, with one positive puncture at $d$, $k$ positive punctures at $a$, and negative punctures at \textbf{b}. The index of a curve $u \in M(d,a^k;\textbf{b})$ is 
\[ind(u) = |d| + k|a| - |\textbf{b}| + k.\]

\noindent
A contractible crossing $a$ is simple if $ind(u) \geq k$ for $u \in M(d,a^k;\textbf{b})$ and $k > 1$. The crossings $b_{i,j}$ of a positive braid are all simple crossings.
\\
\\
In \cite{EHK} it is shown for a simple crossing that the DGA map $\Psi : \mathcal{A}(\Lambda_+) \to \mathcal{A}(\Lambda_-)$ is given by $\Psi(d) = \psi_0(d) + \psi_1(d)$ where:
\begin{itemize}
\item $\psi_0(a) = 1$,
\item $\psi_0(d) = d$ for all $d \in \mathcal{C}(\Lambda_-)$,
\item $\psi_1(a) = 0$, 
\item $\psi_1(d) = \sum_{dim(M(d,a;\textbf{b})) = 1} |M(d,a;\textbf{b})/\mathbb{R}| \cdot \psi_0(\textbf{b})$
for all $d \in \mathcal{C}(\Lambda_-)$.
\end{itemize}
Where $C(\Lambda)$ denotes the set of crossings of $\Lambda$. We note that $\Psi$ is a surjective map.
\\
\\
As an example consider the $(4,3)$ torus braid from Figure $\ref{fig:Torus2}$ whose Legendrian closure is shown in Figure $\ref{fig:plat}$. If we resolve at $b_{1,1}$ then we see that the following crossings have non-zero $\psi_1$:
\begin{itemize}
	\item $\psi_1(b_{2,2}) = b_{1,2}$
	\item $\psi_1(b_{2,3}) = 1$
	\item $\psi_1(c_{4,1}) = c_{4,2}$
	\item $\psi_1(c_{3,1}) = c_{3,2}$
	\item $\psi_1(c_{2,1}) = 1$
\end{itemize}
\subsection{Augmentation Categories.}
We now give a basic overview of $Aug_+(\Lambda)$. We note that everything in Sections 2.3 and 2.4 also holds for $Aug_-(\Lambda)$ with small modifications. For full details the reader should consult \cite{NRSSZ} and \cite{BC}. Let $\Lambda$ be a the Legendrian closure of a positive braid and let $(A,\partial)$ be the Chekanov-Eliashberg DGA associated to $\Lambda$.

\bigskip
\noindent
Our goal will be to construct an $A_{\infty}$ category for $(A,\partial)$. To do this we will need the $m$-copy $\Lambda^{m}$ of $\Lambda$ along with a Morse perturbation $f$ of $\Lambda$. The $m$-copy is obtained by $m-1$ small distinct pushoffs of $\Lambda$ in the Reeb direction of $\mathbb{R}^3_{xyz}$. The Morse perturbation insures that the Reeb chords of $\Lambda^m$ are non-degenerate. We will require that $f$ has $k$ maxima and $k$ minima, where $k$ is the number of strands of $\Lambda$. Let $\Lambda^{m}_f$ denote the $m$-copy together with the Morse perturbation $f$. Label the copies $\Lambda^1,\dots,\Lambda^m$ from top to bottom. Then there is a DGA $(A^m,\partial^m)$ which is generated by the following:
\\
\begin{itemize}
\item Crossings $a^{ij}$ where $1 \leq i,j \leq m$ and $a$ is a crossing of $\Lambda$.
\item Crossings $x^{ij}_l$ where $1 \leq i < j \leq m$ and $1 \leq l \leq k$. These are the crossings associated to the maxima of $f$.
\item Crossings $y^{ij}_l$ where $1 \leq i < j \leq m$ and $1 \leq l \leq k$. These are the crossings associated to the minima of $f$.
\end{itemize}
Throughout the rest of the paper, upper indices will denote components of the $m$-copy and lower indices will be used to distinguish crossings of the knot.
\\
\\
A diagonal augmentation $\boldsymbol{\epsilon}$ of $A^m$ is an $m$-tuple $(\epsilon_1,...,\epsilon_m)$ of augmentations of $\Lambda$ such that $\epsilon_i(b^{jk}) = \delta_{i,j}\delta_{i,k}\epsilon_i(b)$. In this paper all augmentations of $A^m$ will be diagonal.
\\
\\
Inside $A^m$ there are special elements known as composable words. A word $b_{i_1j_1}...b_{i_nj_n}$ is called \textbf{composable} if $j_k = i_{k+1}$ for all $i,j$. We will not describe the differential here; we only need to know that $\partial(b_{ij})$ is a sum of composable words from $i$ to $j$.
\\
\\
We are interested in a linearized version of the above chain complexes. Let $\phi_\epsilon : A_\epsilon^m \to A_\epsilon^m$ be the map defined by $\phi_\epsilon(a) = a + \epsilon(a)$. We define a new differential on $A^m$ by $\partial^m_\epsilon = \phi^{-1}_\epsilon \circ \partial^m \circ \phi_\epsilon$. We note that $\partial^m_\epsilon$ has no constant terms.
\\
\\
Let $C^{ij}$ denote the $\mathbb{Z}/2$ submodule of $(A^m,\partial^m)$ generated by the crossings $a^{ij},x^{ij},y^{ij}$. Then since $\partial$ is made up of composable words we have that $\partial_\epsilon^m|_{C^{ij}}$ breaks into a direct sum of maps of the form 
\[\partial^m_\epsilon : C^{ij} \to C^{ij_1} \otimes ... \otimes C^{i_mj}\]
In particular we are interested in $\partial^m_\epsilon|_{C^{1m}}$. The dual of this map will be denoted $\mu_\epsilon^m$. The multiplication maps $\mu_{\epsilon_1,\epsilon_2}^2$ satisfy $(\mu_{\epsilon_1,\epsilon_2}^2)^2 = 0$, hence $(C^{*12}, \mu_{\epsilon_1,\epsilon_2}^2)$ is a cochain complex where $C^{*12}$ denotes the dual complex.
\\
\\
We can now define the $A_\infty$ augmentation category $Aug_+(\Lambda)$ as follows:
\begin{itemize}
\item \textit{Objects:} The objects are graded augmentations $\epsilon: \Lambda \to \mathbb{Z}_2$.
\item \textit{Morphisms:} The morphisms between two augmentations $\epsilon_1$ and $\epsilon_2$ are the cochain complexes $(C^{*12}, \mu_{\epsilon_1,\epsilon_2}^2)$.
\item \textit{Multiplication Maps:} The multiplication maps are given by the $\mu_\epsilon^m$.
\end{itemize}
\bigskip
As the morphisms of $Aug_+(\Lambda)$ are cochain complexes we can form the cohomology category $HAug_+(\Lambda)$ as follows:
\begin{itemize}
\item \textit{Objects:} The objects are graded augmentations $\epsilon: \Lambda \to \mathbb{Z}_2$
\item \textit{Morphisms:} The morphisms between $\epsilon_1$ and $\epsilon_2$ are the cohomology groups of the chain complex $(C^{*12},\mu_{\epsilon_1,\epsilon_2}^2)$.
\end{itemize}
The cohomology groups for $Aug_-(\Lambda)$ will be denoted $LCH_{\epsilon_1,\epsilon_2}(\Lambda)$
\subsection{Functors.}

We now show how to construct a functor between two augmentation categories. This brings us to the notion of a consistent sequence. A sequence of DGA maps $\Psi^m: (A^m,\partial^m) \to (B^m,\partial^m)$ is called \textbf{consistent} if $\Psi^m(a^{ij})$ is a sum of composable words from $i$ to $j$ in $B^m$. We note that in general, stronger conditions are required for a sequence to be consistent. The reader should consult [NRSSZ] for details.
\\
\\
Given a consistent sequence as above, we construct an $A_\infty$ functor \[F: Aug_+(B) \to Aug_+(A)\] as follows:
\\
\\
On objects we have $F(\epsilon) = \epsilon \circ \Psi$. We denote the submodules of $(A^m,\partial^m)$ and $(B^m,\partial^m)$ as $C^{ij}_A$ and $C^{ij}_B$ respectively.
\\
\\
For each $k$ we need to define maps 
\[F_k: C_B^{*k,k+1} \otimes ... \otimes C_B^{*12} \to C^{*1,k+1}_A.\] 
Consider the diagonal augmentation $\epsilon = (\epsilon_1, ...,\epsilon_{k+1})$ of $B^{k+1}$ and let \[\Psi^{k+1}_\epsilon = \Phi_\epsilon \circ \Psi^{k+1} \circ \Phi^{-1}_{\Psi^*(\epsilon)}.\]
Note that this map contains no constant term. We then define $F_k$ by dualizing the component of $\Psi^{k+1}_\epsilon$ that maps 
\[C_A^{1,k+1} \to C_B^{12} \otimes ... \otimes C_B^{k,k+1}.\]
The $A_\infty$ functor $F$ is defined as the collection of maps $\{F_k\}$.
\\
\\
We note that $F_1$ descends to a well-defined map on cohomology, hence a consistent sequence also induces a functor on cohomology categories.

\section{Results}
\label{Sec:res}
\noindent
Recall that 0-resolution induces a DGA map $\Psi: \mathcal{A}(\Lambda_{+}) \to \mathcal{A}(\Lambda_{-})$. The following lemma gives the structure $\Psi$ for Legendrian closures of positive braids when a crossing in the braid portion is resolved.

\begin{lemma}
	Let $b_{m,n}$ be the crossing to be resolved. Then we claim the following:
	\begin{enumerate}
	\item Let $(i,j) < (m,n)$ in the lexicographic ordering. Then $\psi_1(b_{i,j})$ is a polynomial, possibly with constant term $1$, whose monomials are products of elements of the form $b_{l,k}$ for $(i,j) < (l,k) < (m,n)$.

	\item  Let $(i,j) > (m,n)$ in the lexicographic ordering. Then $\psi_1(b_{i,j})$ is a polynomial, possibly with constant term $1$, whose monomials are products of elements of the form $b_{l,k}$ for $(i,j) > (l,k) > (m,n)$.
	
	\item  $\psi_1(c_{i,j})$ is a polynomial, possibly with constant term $1$, whose monomials are products of elements of the form $b_{l,k}$ for $(l,k) < (m,n)$ and $c_{i,w}$ for $w > j$.
	
	\item  $\psi_1(s_{i,j})$ is a polynomial, possibly with constant term $1$, whose monomials are products of elements of the form $b_{l,k}$ for $(l,k) > (m,n)$ and $a_{i,w}$ for $w < j$.
	\end{enumerate}
\end{lemma}
\noindent
The essential content of the lemma is that the disks which contribute to $\psi_1(d)$ for a crossing $d$ of the link only have negative corners at crossings which are between $d$ and $b_{m,n}$ in the ordering given in the proof of Theorem 3.2.
\\
\\
\textit{Proof:} We prove Case 3, which is slightly more difficult than Cases 1 and 2 and analogous to Case 4. Let $u$ be a disk which contributes to $\psi(c_{i,j})$. We first make the observation that the positive puncture of $u$ at $c_{i,j}$ lies in the bottom left quadrant of $c_{i,j}$. 
\\
\\
Starting at $c_{i,j}$ the top strand of $u$ tavels left, either turning down at some $c$ crossing or none and entering the braid. The bottom strand travels down until it enters the braid. From here both strands continue to the right, possible turning at negative punctures but always continuing to the right. We want to show that the strands cannot extend farther right than $b_{m,n}$. Assume they extend farther to the right. If the strands meet somewhere in the positive braid past $b_{m,n}$ then this is a third positive puncture, which contradicts $u$ contributing to $\psi_1(c_{i,j})$. If they do not meet then they must exit the positive braid on the right side. We then claim that $u$ has a third positive corner at some $s_{i,j}$. Let the bottom strand exit the braid at the $i^{th}$ level. This strand must turn right on or before $s_{i+1,i}$. Then we see that $u$ has a positive corner at $s_{j,j}$ for some $j > i$. \qed
\\
\\
We are now able to prove our main theorem:
\begin{thm}
	Let $\Lambda_+$ be the Legendrian closure of a positive braid and let $\Lambda_-$ be the link obtained after resolving a crossing of $\Lambda_+$. Then there are cohomologically faithful $A_\infty$ \ functors $F^+: Aug_+(\Lambda_-) \to Aug_+(\Lambda_+)$ and $F^-: Aug_-(\Lambda_-) \to Aug_-(\Lambda_+)$ which are induced by 0-resolution. 
\end{thm}
\noindent
\textit{Proof:} We prove Theorem 3.2 for $Aug_+$. The proof for $Aug_-$ is similar. 
\\
\\

Let $\Psi : \mathcal{A}(\Lambda_+) \to \mathcal{A}(\Lambda_-)$ be the map induced by 0-resolution. We will construct a consistent sequence of DGA maps $\Psi^m: A^m(\Lambda_+) \to A^m(\Lambda_-)$. Let $a \in A(\Lambda_+)$. Let $I^{ij}_m$ be the set of all composable sequences from $i$ to $j$ whose superscripts are at most $m$. For a word $b_1b_2...b_l$ in $\Psi(a)$ and $J \in I_m^{ij}$ we denote by $b_1b_2...b_l^J$ the element of $A^m(\Lambda_-)$ corresponding to $b_1b_2...b_l$ indexed by the composable sequence $J$. Then we set
\[\Psi^m(a^{ij}) = \sum_{J \in I_m^{ij}} \Psi(a)^J. \] 
In other words, $\Psi^m(a^{ij})$ is the sum of all possible composable sequences from $i$ to $j$ of $\Psi(a)$.
We also set $\Psi^m$ to be the identity on the mixed chords $x^{ij},y^{ij}$.
\\
\\
This consistent sequence gives rise to an $A_\infty$ functor $F^+$ as in Section 2.4.
\\
\\
To show that $F^+$ is cohomologically faithful, we will show that $F^+_1: C_B^{*12} \to C_A^{*12}$ is injective. We will abuse notation and use $a^{ij}$ to denote both the element of $C^{ij}$ as well as its dual. It is clear the $F^+_1$ is the identity on $x^{12},y^{12}$. We will order the non-Morse elements of $C^{*12}$ in the following way:
\begin{enumerate}
\item $c_{i,j} < c_{m,n}$ if $n < j$, or if $n = j$, $m < i$.

\item $b_{i,j} < b_{m,n}$ if $(i,j) < (m,n)$ in the lexicographic ordering.

\item $s_{i,j} < s_{m,n}$ if $n < j$ or if $n = j$, $i < m$.

\item $c_{i,j} < b_{i,j} < s_{i,j}$ for all $i,j$.
\end{enumerate}
Then by Lemma 3.1, the linear part of $\Psi^{2}_\epsilon$ on the non-Morse chords of $C^{12}_A$ is a block triangular matrix of the form:
\[ \left( \begin{array}{ccccccc}
	1 & * & * & 0 & 0 & 0 & 0 \\
	0 & 1 & * & 0 & 0 & 0 & 0 \\
	0 & 0 & 1 & 0 & 0 & 0 & 0 \\
	0 & 0 & 0 & 0 & 0 & 0 & 0 \\ 
	0 & 0 & 0 & 0 & 1 & 0 & 0  \\
	0 & 0 & 0 & 0 & * & 1 & 0  \\
	0 & 0 & 0 & 0 & * & * & 1  \\
\end{array} \right) \]
\\
The switch from upper triangular to lower triangular occurs at $b_{m,n}$, the resolved crossing. Then $F_1$ is obtained by taking the transpose of this matrix and removing the row and column associated to $b_{m,n}$. This matrix is block triangular and hence is an invertible matrix. Therefore $F_1$ descends to an injective map on cohomology. The proof of Corollary 3.3 will show that this implies F is a cohomologically faithful functor. \qed
\\
\\
Using Theorem 3.2 we can prove:
\begin{cor}
	Let $F^-$ be the functor from Theorem 3.2. Let $\epsilon_1$ and $\epsilon_2$ be augmentations of $\Lambda_-$. Then
	\[LCH^*_{F^-(\epsilon_1),F^-(\epsilon_2)}(\Lambda_+) \simeq LCH^*_{\epsilon_1,\epsilon_2}(\Lambda_-) \oplus \mathbb{Z}/2[0]. \]
\end{cor}
\noindent
\textit{Proof:} Let $C_-$ denote the chain complex associated to $\Lambda_-$ and $C_+$ the chain complex associated to $\Lambda_+$. From the proof of Theorem 3.2 we know that there is an injective map $F_1^-$ from $C_-^*$ to $C_+^*$. This gives rise to a short exact sequence of chain complexes:
\[0 \to C_-^* \xrightarrow{F^-_1} C_+^* \rightarrow C_+^*/C_-^* \to 0.\]
\\
The complex $C^*_+/C^*_-$ consists of one element in degree 0, the crossing $b_{m,n}$ which was resolved. The differential of this complex is identically 0. Therefore we have that its cohomology is isomorphic to $\mathbb{Z}/2$ in degree 0 and is 0 elsewhere. This short exact sequence gives rise to the following long exact sequence in cohomology:
\begin{equation*} 
\begin{aligned}
& 0 \to H_0(C^*_+/C^*_-) \to LCH_0^{F(\epsilon_1),F(\epsilon_2)}(C^*_+)  \to LCH_0^{\epsilon_1,\epsilon_2}(C^*_-) \\
& \to H_1(C^*_+/C^*_-)  \to LCH_1^{F(\epsilon_1),F(\epsilon_2)}(C^*_+) \to LCH_1^{\epsilon_1,\epsilon_2}(C^*_-) \to 0.\\
\end{aligned}
\end{equation*}
\\
\\
It is then immediate from the above sequence that
\[ LCH_0^{F(\epsilon_1),F(\epsilon_2)}(C^*_+) \simeq LCH_0^{\epsilon_1,\epsilon_2}(C^*_-) \oplus \mathbb{Z}/2 .  \]
and
\[ LCH_1^{F(\epsilon_1),F(\epsilon_2)}(C^*_+)   \simeq LCH_1^{\epsilon_1,\epsilon_2}(C^*_-). \] \qed

\subsection{Example}
In this section we show an small example to illustrate the maps for $Aug_-$. Let $\Lambda_+$ be the positive braid with 3 strands whose crossings are labeled $b_{2,1},b_{1,2},b_{2,2},b_{3,2}$. We will resolve $b_{2,1}$ which will result in the $(3,2)$ torus knot union a copy of the unknot. The Legendrian closure of $\Lambda_+$ also has crossings $c_{2,1},c_{3,1},c_{3,2}$ and $s_{2,1},s_{3,1},s_{3,2},s_{1,1},s_{2,2},s_{3,3}$.
\\
\\
$\psi_0$ is the identity on all chords except for $b_{2,1}$ for which $\psi_0(b_{2,1}) = 1$.
The map $\psi_1$ is given as follows:
\begin{itemize}
	\item $\psi_1(b_{2,1}) = 0$
	\item $\psi_1(b_{1,2}) = 0$
	\item $\psi_1(b_{2,2}) = 0$
	\item $\psi_1(b_{3,1}) = 0$
	\item $\psi_1(c_{2,1}) = b_{1,2}$
	\item $\psi_1(c_{3,1}) = 1 + c_{3,2} \cdot b_{1,2}$
	\item $\psi_1(c_{3,2}) = 0$
	\item $\psi_1(s_{3,1}) = 0$
	\item $\psi_1(s_{3,2}) = s_{3,1} + s_{3,1} \cdot b_{2,2} \cdot b_{3,1}$
	\item $\psi_1(s_{3,3}) = s_{3,1} \cdot b_{2,2}$
	\item $\psi_1(s_{2,1}) = 0$
	\item $\psi_1(s_{2,2}) = s_{2,1} + s_{2,1} \cdot b_{2,2} \cdot b_{3,1}$
	\item $\psi_1(s_{1,1}) = 0$ 
\end{itemize}
The map $\Psi: \Lambda_+ \to \Lambda_-$ is the sum of $\psi_0$ and $\psi_1$.
\\
\\
Let $\boldsymbol{\epsilon} = (\epsilon_1,\epsilon_2)$, where $\epsilon_1,\epsilon_2$ are augmentations of $\Lambda_-$. In order to compute $F^-_1$ we will need to know $\Psi^ 2_\epsilon$. $(\psi_0)_\epsilon^2$ is the identity on every crossing except $b_{2,1}$ for which $(\psi_1)_\epsilon^2(b_{2,1}) = 0$.
\\
\\
The non-trivial parts of $(\psi_1)_\epsilon^2$ are as follows:
\begin{itemize}
	\item $c_{2,1} \mapsto b_{1,2}$
	\item $c_{3,1} \mapsto \epsilon_1(c_{3,2}) \cdot b_{1,2} + c_{3,2} \cdot \epsilon_2(b_{1,2})$
	\item $s_{3,2} \mapsto s_{3,1} + \epsilon_1(s_{3,1}) \cdot \epsilon_1(b_{2,2}) \cdot b_{3,1} + \epsilon_1(s_{3,1}) \cdot b_{2,2} \cdot \epsilon_2(b_{3,1}) +
	s_{3,1} \cdot \epsilon_2(b_{2,2}) \cdot \epsilon_2(b_{3,1})$
	\item $s_{3,3} \mapsto \epsilon_1(s_{3,1}) \cdot b_{2,2} + s_{3,1} \cdot \epsilon_2(b_{2,2})$
	\item $s_{2,2} \mapsto \epsilon_1(s_{2,1}) \cdot \epsilon_1(b_{2,2}) \cdot b_{3,1} + \epsilon_1(s_{2,1}) \cdot b_{2,2} \cdot \epsilon_2(b_{3,1}) +
	s_{2,1} \cdot \epsilon_2(b_{2,2}) \cdot \epsilon_2(b_{3,1})$
\end{itemize}
Using the ordering given in Theorem 3.2, this is a matrix of the form given in the proof of Theorem 3.2. Hence removing the row and column corresponding to $b_{2,1}$ and taking the transpose gives $F^-_1$.
\\
\\
The bilinearized Legendrian contact cohomology of the $(3,2)$ torus knot is calculated in [BC]. The extra copy of the unknot contributes a $\mathbb{Z}/2$ factor in degree $1$. Using Theorem 3.3 we see that if $\epsilon_1$ and $\epsilon_2$ are distinct then:
\[LCH_{F^-(\epsilon_1),F^-(\epsilon_2)}^*(\Lambda_+) \simeq (\mathbb{Z}/2)^2[0] \oplus (\mathbb{Z}/2)^2[1].   \]

\end{document}